\documentclass{amsart}
\usepackage{graphicx}
\usepackage{amssymb}

\newtheorem{thm}{Theorem}[section]
\newtheorem{conj}[thm]{Conjecture}
\newtheorem{quest}[thm]{Question}

\newtheorem{lem}[thm]{Lemma}

\theoremstyle{definition}

\theoremstyle{remark}
\newtheorem{rem}[thm]{Remark}

\begin{document}

\title[Foliation]{A proof by foliation that Lawson's cones are $A_{\Phi}$-minimizing}
\author{Connor Mooney}
\address{Department of Mathematics, UC Irvine}
\email{\tt mooneycr@math.uci.edu}
\author{Yang Yang}
\address{Department of Mathematics, UC Irvine}
\email{\tt y.yang@uci.edu}
\begin{abstract}
We give a proof by foliation that the cones over $\mathbb{S}^k \times \mathbb{S}^l$ minimize parametric elliptic functionals for each $k,\,l \geq 1$. We also analyze the behavior at infinity of the leaves in the foliations. This analysis motivates conjectures related to the existence and growth rates of nonlinear entire solutions to equations of minimal surface type that arise in the study of such functionals.
\end{abstract}
\maketitle

\section{Introduction}
A well-known result in the theory of minimal surfaces is that area-minimizing hypersurfaces in $\mathbb{R}^{n+1}$ are smooth when $n \leq 6$, but can have singularities in higher dimensions. An important tool in the theory is the monotonicity formula, which reduces the regularity problem to establishing the existence or non-existence of singular area-minimizing hypercones. Such cones were ruled out by Simons in the case $n \leq 6$ \cite{S}. On the other hand, Bombieri-De Giorgi-Giusti proved that the cone over $\mathbb{S}^3 \times \mathbb{S}^3$ in $\mathbb{R}^8$ is area-minimizing \cite{BDG}.

In this paper we consider the regularity problem for minimizers of parametric elliptic functionals, which generalize the area functional. These assign to an oriented hypersurface $\Sigma \subset \mathbb{R}^{n+1}$ the value
\begin{equation}\label{PEF}
A_{\Phi}(\Sigma) := \int_{\Sigma} \Phi(\nu)\,dA,
\end{equation}
where $\nu$ is the unit normal to $\Sigma$, and $\Phi$ is a one-homogeneous function on $\mathbb{R}^{n+1}$ that is positive and $C^2$ on $\mathbb{S}^n$, and has uniformly convex sub-level sets. Functionals of the form (\ref{PEF}) have attracted recent attention for their applied and theoretical interest (\cite{DDG}, \cite{DMMN}). In particular, they arise in models of crystal surfaces and in Finsler geometry, and the lack of a monotonicity formula for critical points of (\ref{PEF}) presents interesting technical challenges. Almgren-Schoen-Simon proved that the $(n-2)$-dimensional Hausdorff measure of the singular set for a minimizer of (\ref{PEF}) vanishes \cite{ASS}. In particular, minimizers are smooth in the case $n = 2$. Morgan later showed that the cone over $\mathbb{S}^k \times \mathbb{S}^k$ in $\mathbb{R}^{2k+2}$ minimizes a parametric elliptic functional for each $k \geq 1$, by constructing a calibration \cite{Mor}. Thus, there exist singular minimizers of parametric elliptic functionals when $n \geq 3$. 

The purpose of this paper is to prove that the cones over $\mathbb{S}^k \times \mathbb{S}^l$ minimize parametric elliptic functionals for each $k,\,l \geq 1$, by constructing foliations by minimizers in the spirit of \cite{BDG}. To our knowledge, this is the first application of the foliation approach for integrands other than area (that is, $\Phi|_{\mathbb{S}^n} = 1$).

\begin{rem}
The examples here and in \cite{Mor} show that the best regularity result possible for minimizers of (\ref{PEF}) is that the singular set has e.g. locally bounded $(n-3)$-dimensional Hausdorff measure. It remains an interesting open problem to determine the maximum possible dimension of the singular set for a minimizer of (\ref{PEF}). See e.g. \cite{W}, pg. $686$, for further discussion of this problem.
\end{rem}

The approach of constructing a foliation by minimizers has several advantages. The first is that it removes some of the guesswork involved in constructing a calibration. Indeed, the approach involves solving a nonlinear ODE, which we show is possible provided the integrand $\Phi$ satisfies analytic conditions that are straightforward to check (see Lemma \ref{Technical}). The second is that the behavior at infinity of a leaf in the foliation gives quantitative information that is useful in the study of the closely related Bernstein problem for graphical minimizers. When a critical point of (\ref{PEF}) can be written as the graph of a function $u$, we say that $u$ solves an equation of minimal surface type. An interesting question is:

\begin{quest}\label{BP}
Are entire solutions to equations of minimal surface type in $\mathbb{R}^n$ necessarily linear?
\end{quest} 

\noindent For the area functional, the answer to Question \ref{BP} is ``yes" if $n \leq 7$ (\cite{F}, \cite{DeG}, \cite{A}, \cite{S}) and ``no" if $n \geq 8$ (\cite{BDG}). For general parametric elliptic functionals, the answer is ``yes" when $n \leq 3$ (\cite{J}, \cite{Si2}) and was recently shown to be ``no" when $n \geq 6$ (\cite{M}). Our main theorem sheds light on the remaining open cases $n = 4,\,5$. More precisely, it suggests the existence of nonlinear entire solutions to equations of minimal surface type in dimension $n \geq 4$ that grow sub-quadratically at infinity (see Conjecture \ref{ExistenceConj} and the discussion after its statement).

For future reference we state our main result here. We fix $k,\, l \geq 1$, we let $x \in \mathbb{R}^{k+1}$ and $y \in \mathbb{R}^{l+1}$, and we define the Lawson cones $C_{kl}$ over $\mathbb{S}^k \times \mathbb{S}^l$ by
$$C_{kl} := \{|x| = |y|\} \subset \mathbb{R}^{k+l+2}.$$

\begin{thm}\label{Main}
For each $k,\,l \geq 1$, there exist parametric elliptic functionals $A_{\Phi}$ such that $\Phi$ is analytic away from the origin, and each side of the cone $C_{kl}$ is foliated by analytic minimizers of $A_{\Phi}$. In particular, $C_{kl}$ minimizes $A_{\Phi}$.
\end{thm}

\begin{rem}
The foliation is generated by dilations of a pair of critical points of $A_{\Phi}$, each of which lies on one side of $C_{kl}$ and is asymptotic to $C_{kl}$ at infinity. We discuss the precise asymptotic behavior in Section \ref{Discussion}.
\end{rem}

\begin{rem}
The $A_{\Phi}$-minimality of $C_{kl}$ seems to be new in the cases $k \neq l$ and $k + l \leq 5$, or $k + l = 6$ and $\min\{k,\,l\} = 1$. The cases $k = l \leq 2$ were treated in \cite{Mor}, and the remaining examples ($k + l = 6$ and $\min\{k,\,l\} \geq 2$ or $k + l \geq 7$) are known to minimize area, up to making an affine transformation (\cite{L}, \cite{Sim}).
\end{rem}

\begin{rem}
A proof of minimality by foliation gives rise to a proof by calibration through the following observation: If we denote by $\nu(z)$ the unit normal to the leaf that passes through $z$, then the vector field $\nabla \Phi(\nu)$ is a calibration on $\mathbb{R}^{k+l+2}$. For a discussion of this connection in the area case see e.g. \cite{Dav}.
\end{rem}

The paper is organized as follows. In Section \ref{MainSection} we prove Theorem \ref{Main}. We reduce the problem to the careful analysis of a certain nonlinear second-order ODE, using the symmetries of $C_{kl}$. More precisely, we choose integrands that (like $C_{kl}$) are invariant under rotations in $x$ and $y$, and
we search for critical points $\Sigma_{kl}$ of (\ref{PEF}) that share these symmetries and can be written in the form $\{|y| = \sigma(|x|)\}$ or $\{|x| = \sigma(|y|)\}$ for some function $\sigma$ of one variable. The condition that $\Sigma_{kl}$ is a critical point gives rise to a nonlinear second-order ODE for $\sigma$. The heart of our construction is Lemma \ref{Technical}, which gives conditions on the integrand that guarantee the existence of solutions to the ODE with the desired properties. In particular, after a change of variable we can view the ODE as a nonlinear first-order autonomous system. The conditions we impose on the integrand guarantee that the solution trajectory is trapped in a region of the plane that corresponds to a function $\sigma$ that defines a foliation leaf $\Sigma_{kl}$ which is asymptotic to $C_{kl}$, and approaches $C_{kl}$ at a precise rate. We remark that our approach quickly recovers the foliation by area-minimizing hypersurfaces of each side of the Simons cone $C_{kk}$ when $k \geq 3$ (see Remark \ref{AreaCase}). Finally, in Section \ref{Discussion} we discuss the behavior at infinity of the leaves in the foliation given by Theorem \ref{Main}, and the implications for Question \ref{BP}. In particular, we state conjectures concerning the existence and growth rates of nonlinear global solutions to equations of minimal surface type in $\mathbb{R}^{k+l+2}$ for each $k,\,l \geq 1$, and we compare these conjectures with what is known about the minimal surface equation.


\section*{Acknowledgements}
We are grateful to the anonymous referee for helpful comments which improved the exposition. This research was supported by NSF grant DMS-1854788.

\section{Proof of Theorem \ref{Main}}\label{MainSection}

\subsection{Integrand Notation}
We choose integrands $\Phi$ that depend only on $|x|$ and $|y|$. We define them by a pair of one-variable functions $\phi$ and $\psi$ as follows:
\begin{equation}\label{PhiNotation}
\Phi(x,\,y) = \begin{cases}
|y|\phi\left(\frac{|x|}{|y|}\right), \quad |y| \geq |x| \\
|x|\psi\left(\frac{|y|}{|x|}\right), \quad |x| > |y|.
\end{cases}
\end{equation}
The functions $\phi$ and  $\psi$ will be chosen to be positive, smooth, even, and locally uniformly convex on $\mathbb{R}$.

\subsection{Foliation Leaf Notation}
\begin{figure}
 \begin{center}
    \includegraphics[scale=0.6, trim={35mm 30mm 10mm 0mm}, clip]{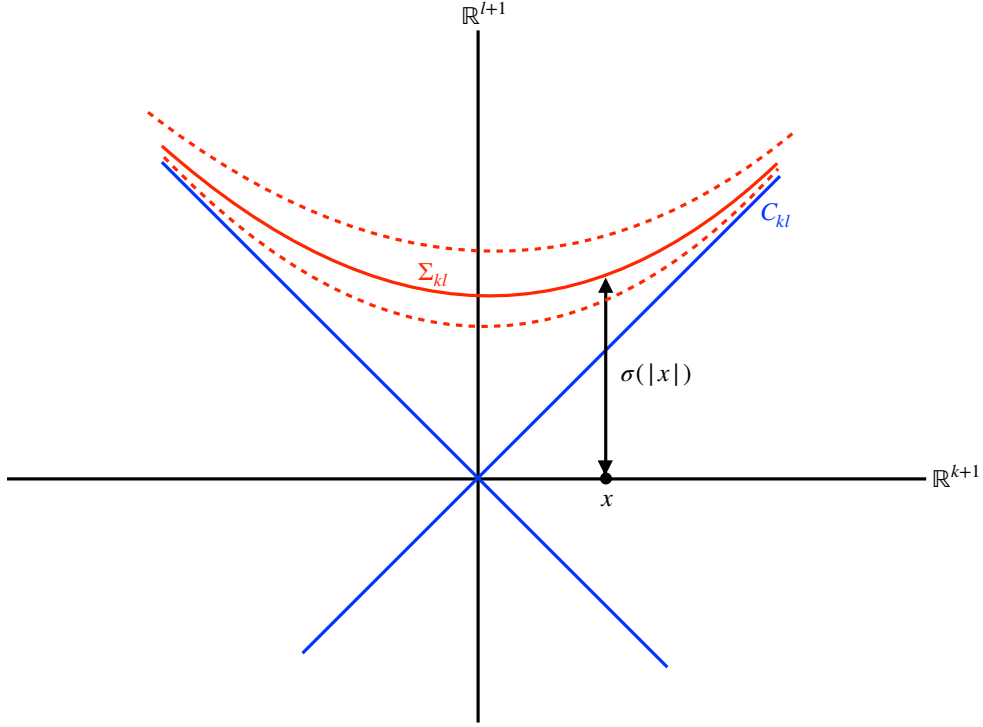}
\caption{The dilations of $\Sigma_{kl}$ foliate one side of $C_{kl}$.}
\label{LeafPic}
\end{center}
\end{figure}

Having fixed an appropriate choice of $\Phi$, we will show that there exists a critical point of $A_{\Phi}$ of the form
\begin{equation}\label{Revolution}
\Sigma_{kl} = \{|y| = \sigma(|x|)\} \subset \{|y| > |x|\},
\end{equation}
where $\sigma$ is smooth, even, convex, asymptotic to $|.|$, and $|\sigma'| < 1.$ The dilations of $\Sigma_{kl}$ are then minimizers of $A_{\Phi}$, and they foliate one side of $C_{kl}$ (namely $\{|y| > |x|\}$), see Figure \ref{LeafPic}. A similar procedure will give a foliation of the other side $\{|x| > |y|\}$ by minimizers of $A_{\Phi}$.

\subsection{Euler-Lagrange ODE}
For hypersurfaces of the form (\ref{Revolution}) and integrands of the form (\ref{PhiNotation}), the condition that $\Sigma_{kl}$ is
a critical point of $A_{\Phi}$ is equivalent to the nonlinear second-order ODE
\begin{equation}\label{EL}
\sigma''(t) + kP(\sigma'(t))\frac{\sigma'(t)}{t} + lQ(\sigma'(t))\frac{1}{\sigma(t)} = 0,
\end{equation}
where
\begin{equation}\label{PQ}
P(s) := \frac{\phi'(s)}{s\phi''(s)}, \quad Q(s) := \frac{s\phi'(s) - \phi(s)}{\phi''(s)}.
\end{equation}
This follows from the first variation formula 
$$tr(D^2\Phi(\nu(z))II(z)) = 0$$ 
for critical points of (\ref{PEF}), where $II$ denotes the second fundamental form of $\Sigma$ and $\nu$ denotes the unit normal. One can also use the symmetries of $\Sigma_{kl}$ and $\Phi$ to reduce the problem to taking the first variation of the one-variable integral
$$A_{\Phi}(\Sigma_{kl}) = const.\int t^k\sigma^l(t)\phi(\sigma'(t))\,dt.$$

In the following technical lemma we show that there exists a global solution to (\ref{EL}) with the desired properties, provided $\phi$ satisfies certain analytic conditions. We will later give examples of $\phi$ that satisfy these conditions. To state the lemma we define for a smooth function $\varphi$ on $\mathbb{R}$ the function $E_{kl}(\varphi)$ by
\begin{equation}\label{TrappingQuantity}
\begin{split}
 E_{kl}(\varphi)(s) := l\frac{k+l-1}{k+l+1}\varphi(s) - \left(k + \left(l - 2\frac{k+l}{k+l+1}\right)s\right)\varphi'(s) \\
 - \left(\frac{k+l+1}{2}-s\right)(1-s)\varphi''(s).
 \end{split}
\end{equation}

\begin{lem}\label{Technical}
Assume that $\phi(s)$ is a smooth, even, uniformly convex function on $\mathbb{R}$ that satisfies
\begin{equation}\label{1Jet}
\phi(1) = 1, \quad \phi'(1) = \frac{l}{k+l},
\end{equation}
and in addition that
\begin{equation}\label{Trapping}
E_{kl}(\phi)(s) \geq \kappa(1-s)
\end{equation}
for some $\kappa > 0$ and all $s \in [0,\,1]$. Then there exists a global smooth, even, convex solution $\sigma$ to the ODE (\ref{EL}) that satisfies the initial conditions 
\begin{equation}\label{ICs}
\sigma(0) = 1, \quad \sigma'(0) = 0
\end{equation} 
and in addition satisfies $\sigma(t) > |t|,\, |\sigma'(t)| < 1$ for all $t$, and
\begin{equation}\label{Asymptotics}
\sigma(t) = |t| + a|t|^{-\mu} + o(|t|^{-\mu})
\end{equation}
as $|t| \rightarrow \infty$ for some $a > 0$, where
\begin{equation}\label{Rate}
\mu := \frac{k+l-1}{2} - \sqrt{\left(\frac{k+l-1}{2}\right)^2 - \frac{kl}{\phi''(1)(k+l)}}.
\end{equation}
\end{lem}

\begin{rem}\label{QuantTrapping}
It is straightforward to check that any function $\phi$ satisfying the conditions (\ref{1Jet}) and (\ref{Trapping}) automatically satisfies the inequality
\begin{equation}\label{SecondDeriv}
\phi''(1)  - \frac{4kl}{(k+l)(k+l-1)^2} > 0,
\end{equation}
so $\mu$ is well-defined. Conversely, any choice of $\phi$ that satisfies (\ref{1Jet}) and (\ref{SecondDeriv}) also satisfies (\ref{Trapping}) for $s \in [1-\delta,\, 1]$, where $\kappa > 0$ and $\delta > 0$ depend only on the left side of (\ref{SecondDeriv}) and 
$\|\phi\|_{C^3([-1,\,1])}$.
\end{rem}

\begin{proof}[{\bf Proof of Lemma \ref{Technical}}]
Standard ODE theory gives the short-time existence of a solution to (\ref{EL}) with the desired properties in a neighborhood of $0$ (see Remark \ref{STE} below). To proceed we rewrite (\ref{EL}) 
as an autonomous first-order system. In terms of the quantities
\begin{equation}\label{NewVars}
w(\tau) := e^{-\tau}\sigma(e^{\tau}), \quad z(\tau) := \sigma'(e^{\tau}),
\end{equation}
the second-order ODE (\ref{EL}) becomes:
\begin{equation}\label{NonlinearSystem}
\left(\begin{array}{c}
w'\\
z'\\
\end{array} \right) = \left(\begin{array}{c}
-w + z\\
-l\frac{Q(z)}{w} - kzP(z)\\
\end{array} \right) := {\bf V}(w,\,z).
\end{equation}
We denote the components of the vector field ${\bf V}$ by $V^i,\, i = 1,\,2$, and the solution curve $(w(\tau),\,z(\tau))$ by $\Gamma(\tau)$. The only zero of ${\bf V}$ in the closure of the infinite half-strip 
$$\Omega := \{w > 1\} \cap \{0 < z < 1\}$$ 
occurs at $(1,\,1)$ (here we used (\ref{1Jet})). In addition, the linearization of (\ref{NonlinearSystem}) around the zero $(1,\,1)$ has the form $X' = M\,X$, where
\begin{equation}
M = \begin{pmatrix}
-1 & 1 \\
-\frac{kl}{\phi''(1)(k+l)} & -k-l
\end{pmatrix}.
\end{equation}
The eigenvalues of $M$ are
\begin{equation}\label{Eigvals}
\lambda_{\pm} = -\frac{k+l+1}{2} \pm \sqrt{\left(\frac{k+l-1}{2}\right)^2 - \frac{kl}{\phi''(1)(k+l)}},
\end{equation}
and these eigenvalues correspond to directions with slopes $1 + \lambda_{\pm}$. 

We claim that $\Gamma$ is contained in the region $R \subset \Omega$ bounded by the curves 
\begin{align*}
\Gamma_1 := &\{z = 0\}, \quad \Gamma_2 := \{V^2 = 0\} = \left\{w = \frac{l}{k}\left(\frac{\phi}{\phi'} - z\right)\right\}, \text{ and } \\
&\Gamma_3 := \left\{(z -1) = \left(1 + \frac{\lambda_+ + \lambda_-}{2}\right)(w-1)\right\}.
\end{align*}
We first note that $V^2 > 0$ on $\Gamma_1 \cap \{w > 0\}$ using that $Q(0) < 0$. Next, by the uniform convexity of $\phi$, the curve $\Gamma_2 \cap \Omega$ is a graph over the positive $z$-axis with negative slope, and furthermore $V^1 < 0$ in $\Omega$ (in particular, on $\Gamma_2 \cap \Omega$). Finally, after a calculation using the definitions (\ref{PQ}) of $P$ and $Q$ and inequality (\ref{Trapping}), we have 
$$\frac{V^2}{V^1} < 1 + \frac{\lambda_+ + \lambda_-}{2}$$
on $\Gamma_3 \cap \Omega$. Indeed, the preceding inequality can be written
$$-lQ(z) - kwzP(z) > \frac{k+l-1}{2}w(w-z)$$
on $\Gamma_3 \cap \Omega$. Using that
$$\Gamma_3 \cap \Omega = \left\{w = \frac{k+l+1}{k+l-1} - \frac{2}{k+l-1}z,\, z \in (0,\,1)\right\}$$
and that
$$-Q(z) = \frac{\phi(z)-z\phi'(z)}{\phi''(z)}, \quad -zP(z) = -\frac{\phi'(z)}{\phi''(z)}$$
the previous inequality becomes
\begin{align*}
l(\phi-z\phi') - &k\left(\frac{k+l+1}{k+l-1} - \frac{2}{k+l-1}z\right)\phi' \\
&> \frac{k+l-1}{2}\left(\frac{k+l+1}{k+l-1} - \frac{2}{k+l-1}z\right)\frac{k+l+1}{k+l-1}(1-z)\phi'' \\
&= \frac{k+1+1}{k+l-1}\left(\frac{k+l+1}{2} - z\right)(1-z)\phi''.
\end{align*} 
After regrouping terms, we see that this inequality holds for $z \in (0,\,1)$ using (\ref{Trapping}). We conclude that $\Gamma_2$ and $\Gamma_3$ meet only at $(1,\,1)$, and
${\bf V}$ points towards the interior of $R$ on each of the curves $\{\Gamma_i\}_{i = 1}^3$ (see Figure \ref{TrappingPic}). It only remains to argue that $\Gamma(\tau) \in R$ for all $\tau << 0$, which holds by the initial convexity of $\sigma$ (note that $V^2 < 0$ ``above" $\Gamma_2$ in $\Omega$).

\begin{figure}
 \centering
    \includegraphics[scale=0.55, trim={40mm 40mm 30mm 20mm}, clip]{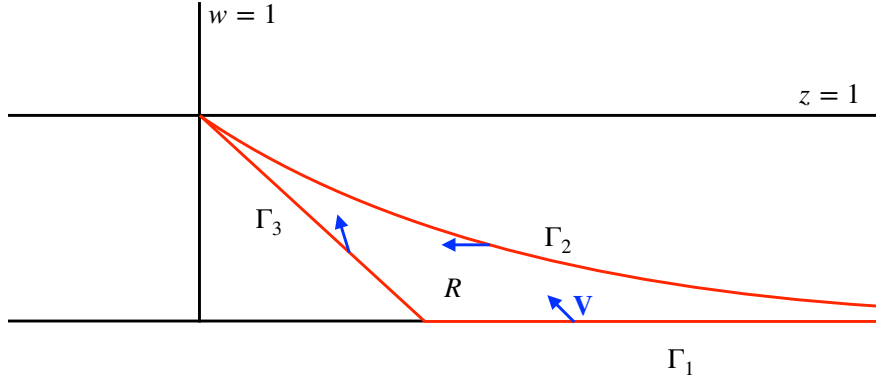}
\caption{The solution curve is contained in the region $R$ bounded by $\Gamma_1,\, \Gamma_2$ and $\Gamma_3$.}
\label{TrappingPic}
\end{figure}

Since $R \subset \Omega \cap \{V^2 > 0\},$ we conclude that $\sigma > |.|,\, |\sigma'| < 1$ and that $\sigma$ is convex. The asymptotic behavior (\ref{Asymptotics}) follows from the linear analysis. Indeed, the region $R$ excludes the line with slope $1+\lambda_{-}$ that goes through $(1,\,1)$. Thus, in the expansion of $\Gamma$ for $\tau$ large, the coefficient of the principal eigenvector of the linearized operator at $(1,\,1)$ (which corresponds to the eigenvalue $\lambda_{+}$ and has slope $-\mu = 1+\lambda_+$) is nonzero, completing the proof.
\end{proof}

\begin{rem}\label{STE}
We could not find a precise reference for short-time existence, so for completeness we sketch the argument. We first rewrite (\ref{EL}) in divergence form:
$$[t^k\sigma^l\phi'(\sigma')]' = lt^k\sigma^{l-1}\phi(\sigma').$$
We are thus looking for a continuous function $\sigma'$ on an interval $[0,\,t_0]$ such that $\sigma'(0) = 0$ and 
\begin{align*}
\sigma'(t) &= (\phi^*)'\left(\frac{l}{t^k\left(1+\int_0^t \sigma'(s)\,ds\right)^l}\int_0^t s^k\left(1+\int_0^s \sigma'(\tau)\,d\tau\right)^{l-1}\phi(\sigma'(s))\,ds\right) \\
&:= G(\sigma')(t),
\end{align*}
where $\phi^*$ is the Legendre transform of $\phi$. For $t_0 > 0$ small, the operator $G$ is a contraction mapping on the space of continuous functions on $[0,\,t_0]$ that vanish at $0$ and are bounded by $1$ in the $C^0$ norm. A fixed point argument then gives the existence of a function $\sigma \in C^1[0,\,t_0] \cap C^{\infty}(0,\,t_0)$ that solves (\ref{EL}) on $(0,\,t_0)$ and satisfies $\sigma(0) = 1,\, \sigma'(0) = 0$.
The higher regularity of $\sigma$ follows from the observation that $\Sigma_{kl} = \{|y| = \sigma(|x|)\}$ can be locally written over its tangent planes as a $C^1$ graph that solves an equation of minimal surface type (and is thus smooth, see e.g. \cite{GT}). Finally, the equation (\ref{EL}) itself gives that 
$$\sigma''(0) = \frac{l\phi(0)}{(k+1)\phi''(0)} > 0,$$
so $\sigma$ is convex near $0$, concluding the argument.
\end{rem}

\begin{rem}\label{AreaCase}
In the case of the area functional we have that 
\begin{equation}\label{AreaIntegrand}
\phi(s) = \sqrt{1 + s^2}.
\end{equation}
Our approach recovers the foliation of each side of the cone $C_{kk}$ by area-minimizing hypersurfaces when $k \geq 3$, as follows. First, a short calculation shows that the function $\sigma_0(t) := (1 + t^4)^{1/4}$ is a super-solution of the ODE (\ref{EL}) corresponding to the area integrand (\ref{AreaIntegrand}) when $k \geq 3$. A similar calculation was performed in \cite{DP} to construct a sub-calibration for $C_{kk}$. If we take $\Gamma_2$ as above and we take
$\Gamma_3 = \{(e^{-\tau}\sigma_0(e^{\tau}),\,\sigma_0'(e^{\tau})), \, \tau \in \mathbb{R}\}$, then a similar argument to the one above shows that the solution trajectory to the associated autonomous system is trapped between $\Gamma_2$ and $\Gamma_3$, giving an exact solution $\sigma$ of (\ref{EL}) with the desired properties.
\end{rem}

\subsection{Proof of Main Theorem}

In this final subsection we choose the functions $\phi,\,\psi$ that define $\Phi$, and we apply Lemma \ref{Technical} to prove Theorem \ref{Main}.

\begin{proof}[{\bf Proof of Theorem \ref{Main}}]
We first indicate how to choose integrands $\Phi$ that are $C^{2,\,1}$ away from the origin, defined through the notation (\ref{PhiNotation}). For $p, q > 2$ to be chosen later we take
\begin{equation}\label{PhiPsiChoice}
\begin{split}
\phi(s) = 1-\frac{l}{p(k+l)} + \frac{l}{p(k+l)}|s|^p, \\
\psi(s) = 1-\frac{k}{q(k+l)} + \frac{k}{q(k+l)}|s|^q,
\end{split}
\end{equation}
up to making small perturbations near $s = 0$ so that $\phi$ and $\psi$ are smooth and uniformly convex. It is straightforward to check that if $p$ and $q$ are related by
\begin{equation}\label{Compatibility}
l(p-1) = k(q-1),
\end{equation}
then $\Phi$ is $C^{2,\,1}$ away from the origin. We note that $\phi$ satisfies (\ref{1Jet}), and we will verify that provided $p$ is sufficiently large, then $\phi$ also satisfies the desired inequality (\ref{Trapping}). Away from a small neighborhood of $s = 0$, where we perturbed (\ref{PhiPsiChoice}) and the inequality (\ref{Trapping}) is obvious, the inequality $E_{kl}(\phi)(s) > 0$ becomes
\begin{equation}\label{PhiTrapping}
\begin{split}
\left[(p-1)\left(\frac{k+l+1}{2}-s\right) + ks\right](1-s) \\
< \frac{(k+l-1)(k+l-l/p)}{k+l+1}(s^{2-p}-s^2).
\end{split}
\end{equation}
Denote the left side of (\ref{PhiTrapping}) by $L(s)$ and the right side by $R(s)$. Since $L$ is quadratic in $s$ and $R''$ is decreasing in $s$, it suffices to prove the inequalities
$$R'(1) < L'(1), \quad L'' \leq R''(1).$$
The first inequality holds provided 
\begin{equation}\label{FirstOrder}
p-1 > \frac{4k}{(k+l-1)^2},
\end{equation}
in agreement with Remark (\ref{QuantTrapping}). The second one holds provided
\begin{equation}\label{SecondOrder}
\begin{split}
(p-1)^2 - \left(\frac{k+2l+2}{k+l} + \frac{4}{(k+l)(k+l-1)}\right)(p-1) \\
 + \frac{4k}{(k+l)(k+l-1)} \geq 0.
\end{split}
\end{equation}
Both (\ref{FirstOrder}) and (\ref{SecondOrder}) hold e.g. when $p \geq 6$, regardless of $k,\,l \geq 1$. When $p \geq 6$ the inequality (\ref{SecondDeriv}) also holds, so by Remark (\ref{QuantTrapping}) the desired inequality (\ref{Trapping}) holds for some $\kappa > 0$. 

Up to exchanging $k$ and $l$, the function $\psi$ satisfies (\ref{1Jet}), and a similar analysis shows that $E_{lk}(\psi) \geq \kappa(1-s)$ for some $\kappa > 0$ and all $s \in [0,\,1]$ if $q \geq 6$.
We conclude using Lemma \ref{Technical} that each side of $C_{kl}$ is foliated by smooth critical points of $A_{\Phi}$ when we choose $\phi,\,\psi$ as above with $p,\,q \geq 6$, and furthermore $\Phi \in C^{2,\,1}(\mathbb{S}^{k+l+1})$ provided $p$ and $q$ are chosen such that (\ref{Compatibility}) holds as well.

We now explain how the integrand can be made analytic on $\mathbb{S}^{k+l+1}$, by perturbing the $C^{2,\,1}$ integrand constructed above. We first improve to smooth. Take $\phi$ and $\psi$ as above, and let
$$\tilde{\phi}(s) = s\psi(1/s)$$
for $s \leq 1$. We glue $\phi$ to $\tilde{\phi}$ near $s = 1$ by taking the convex combination 
$$\bar{\phi} := \eta_{\delta}\phi + (1-\eta_{\delta})\tilde{\phi},$$ 
where $\eta_{\delta}$ is a smooth function that transitions from $1$ to $0$ in the interval $[1-2\delta,\, 1-\delta]$ for $\delta > 0$ to be chosen, and satisfies 
$$\|\eta_{\delta}\|_{C^m(\mathbb{R})} \leq C_m\delta^{-m}$$ 
with $C_m$ independent of $\delta$. Since $\tilde{\phi}$ and $\phi$ agree to second order at $s = 1$, the inequality (\ref{Trapping}) holds for $\bar{\phi}$ away from $[1-2\delta,\, 1-\delta]$ provided $\delta$ is small (see Remark (\ref{QuantTrapping})). Furthermore, we have
$$|\tilde{\phi}^{(m)}(s) - \phi^{(m)}(s)| \leq C_m(1-s)^{3-m}$$
for each $m \leq 2$ and $s \in [1/2,\,1]$. It follows that
$$|E_{kl}(\bar{\phi}) - E_{kl}(\phi)| \leq C\delta^2$$
in $[1-2\delta,\,1-\delta]$.
Since $E_{kl}(\phi) \geq \kappa \delta$ in this interval, the inequality (\ref{Trapping}) holds for $\bar{\phi}$ when $\delta$ is small, up to reducing $\kappa$ slightly. After replacing $\phi$ by $\bar{\phi}$ (and keeping $\psi$ the same), we obtain a new integrand that is smooth on $\mathbb{S}^{k+l+1}$ and by Lemma \ref{Technical} satisfies the desired properties.

Finally, we indicate how to improve the regularity from smooth to analytic. We start with a smooth choice of integrand $\Phi$ as constructed above. Using the symmetries of $\Phi$ we may view it as a smooth function on $\mathbb{S}^1$. We approximate this function by the partial sums $S_N$ of its Fourier series with $N$ terms. We add small correctors of the form $a_N + b_N\cos(2\theta) + c_N\cos(4\theta)$ to $S_N$ to obtain new approximations $T_N$, with $a_N,\,b_N,\,c_N$ chosen such that $T_N$ agrees to second order with $\Phi$ at $\theta = \pi/4$. Since $S_N$ converge uniformly in $C^m$ to $\Phi$ for any $m$, the functions $T_N$ do as well. It follows that the one-homogeneous extensions of $T_N$ to $\mathbb{R}^2$ (which we now identify with $T_N$) have uniformly convex sub-level sets for $N$ large. Since $T_N$ agree to second order with $\Phi$ on the diagonals, Remark (\ref{QuantTrapping}) implies that the conditions (\ref{1Jet}) and (\ref{Trapping}) hold for the function obtained by restricting $T_N$ to the horizontal lines tangent to $\mathbb{S}^1$ when $N$ is large. The same holds (with $k$ and $l$ exchanged) for the restriction of $T_N$ to the vertical lines tangent to $\mathbb{S}^1$. Hence, after replacing $\Phi(x,\,y)$ with $T_{N}(|x|,\,|y|)$ for $N$ large, we obtain an integrand that is analytic on $\mathbb{S}^{k+l+1}$ and by Lemma \ref{Technical} satisfies the desired properties. 
\end{proof}

\section{Discussion}\label{Discussion}

In this section we discuss the implications of the analysis in Section \ref{MainSection} for Question \ref{BP}. The discussion is motivated by the examples of entire minimal graphs constructed in \cite{BDG} and \cite{Si1}. Those examples are asymptotic to area-minimizing cones of the form $K \times \mathbb{R}$, where $K$ is the Simons cone in \cite{BDG}, and any one of a large family of area-minimizing cones with isolated singularities in \cite{Si1}. In all cases, each side of $K$ is foliated by smooth area-minimizing hypersurfaces. These are closely related to the level sets of the functions $u$ that define the entire minimal graphs. More precisely, each level set of $u$ is a graph over $K$ outside of some ball, with the same leading-order asymptotic behavior at infinity as a leaf in the foliation. Furthermore, if the distance between a leaf in the foliation and $K$ on $\partial B_r$ behaves like $r^{-\mu}$ as $r \rightarrow \infty$, then $\sup_{B_r}|\nabla u| \sim r^{\mu}$. 

In view of this discussion we conjecture:

\begin{conj}\label{ExistenceConj}
For any integrand $\Phi$ as constructed in Theorem \ref{Main}, there exists an elliptic extension of $\Phi$ to $\mathbb{R}^{k+l+3}$, and a nonlinear global solution to the corresponding equation of minimal surface type in $\mathbb{R}^{k+l+2},$ whose graph is asymptotic to $C_{kl} \times \mathbb{R}$. Moreover, the gradient of this solution grows at the same rate that the leaves in the foliation associated to $\Phi$ approach $C_{kl}$.
\end{conj}

\noindent The proof of Theorem \ref{Main} shows that for any $\mu \in (0,\, \mu_{kl})$, we can choose integrands such that each side of $C_{kl}$ is foliated by minimizers whose distance from $C_{kl}$ on $\partial B_r$ behaves like $r^{-\mu}$, where
\begin{equation}\label{MaxRate}
\mu_{kl} = \frac{k+l-1}{2} - \sqrt{\left(\frac{k+l-1}{2}\right)^2 - \frac{\min\{k,\,l\}}{5}}.
\end{equation}
The formula (\ref{MaxRate}) comes from (\ref{Rate}) and noting that, when choosing $\phi$ and $\psi$, we could take any exponents $p$ and $q$ such that (\ref{Compatibility}) holds and $p,\,q \geq 6$. Thus, Conjecture (\ref{ExistenceConj}) predicts that for any $\mu \in (0,\,\mu_{kl})$, there exist global solutions to equations of minimal surface type in $\mathbb{R}^{k+l+2}$ whose graphs are are asymptotic to $C_{kl} \times \mathbb{R}$, and have maximum gradient in $B_r$ growing like $r^{\mu}$.

\begin{rem}
The first-named author showed in \cite{M} that when $k = l = 2$, the graph of $u = |x|^2 - |y|^2$ (which is asymptotic to $C_{22} \times \mathbb{R}$) minimizes a parametric elliptic functional $A_{\Psi}$, and each level set of $u$ minimizes $A_{\Phi}$, where $\Phi = \Psi|_{\{x_7 = 0\}}$. The perspective in that work is quite different, and the proof is based on solving a linear hyperbolic equation to construct $\Psi$. However, the discussion at the end of \cite{M} shows that this strategy could be challenging to implement when $2 \leq k+l \leq 3$. In these cases, Question \ref{BP} may instead yield to a combination of the approaches from \cite{BDG} and \cite{M}.
\end{rem}

To conclude we discuss the gradient growth rates of the solutions predicted by Conjecture \ref{ExistenceConj}. We first consider the possibility of constructing solutions with fast growth. In the case $k = l = 1$, a closer inspection of inequalities (\ref{FirstOrder}) and (\ref{SecondOrder}) shows that we can take any $p = q > 5$ when defining $\phi$ and $\psi$. This corresponds to the ``optimal" value $\mu_{11} = \frac{1}{2}$. One may hope to show that (\ref{MaxRate}) can be improved to $\mu_{kl} = (k + l - 1)/2$ for arbitrary $k$ and $l$, which corresponds to choices of $\phi$ such that inequality (\ref{SecondDeriv}) tends to equality. However, we suspect that this is not possible. Indeed, if $k=l$ is large, this corresponds to a small value of $\phi''(1)$. It would follow that the integrand $\Phi$ is larger on $\mathbb{S}^k \times \mathbb{S}^k$ than at nearby points on $\sqrt{2}\,\mathbb{S}^{2k+1}$, in which case perturbations of $C_{kk}$ could likely decrease its energy. Since the quantity (\ref{MaxRate}) is bounded above independently of $k$ and $l$, it does not seem likely that the examples from Theorem \ref{Main} can give rise to solutions to equations of minimal surface type with arbitrarily fast gradient growth.

On the other hand, Conjecture \ref{ExistenceConj} predicts the existence of global solutions to equations of minimal surface type with very slow gradient growth, namely 
$$\sup_{B_r}|\nabla u| \sim r^{\mu}$$ 
with $\mu > 0$ small. However, global solutions to equations of minimal surface type with bounded gradient are linear. This is a consequence of the De Giorgi-Nash-Moser theorem. Indeed, if the gradient of a global solution $u$ to an equation of minimal surface type is bounded, then each derivative $u_e$ is a global bounded solution to a uniformly elliptic equation in divergence form (with ellipticity constants depending on $\|\nabla u\|_{L^{\infty}}$). Applying e.g. Theorem $8.22$ from \cite{GT} and sending $R_0 \rightarrow \infty$ (keeping in mind the remark at the end of Section $8.9$), we see that $u_e$ is constant. We thus expect that the ellipticity of the integrands from Conjecture \ref{ExistenceConj} will degenerate as $\mu$ tends to zero, and we conjecture a ``quantitative" version of the rigidity result for solutions with bounded gradient:

\begin{conj}\label{GrowthConj}
Let $u$ be a global solution to an equation of minimal surface type on $\mathbb{R}^n$, corresponding to a functional $A_{\Psi}$. Then for some $\epsilon(n,\,\Psi) > 0$,
$$\sup_{B_r} |\nabla u| = O(r^{\epsilon}) \Rightarrow u \text{ is linear.}$$
\end{conj}

\noindent In \cite{EH} the authors give a beautiful proof of Conjecture \ref{GrowthConj} for the area functional, for any $\epsilon < 1$ and in arbitrary dimension $n$. The proof in \cite{EH} depends on precise constants in the Simons inequality for the Laplacian of the second fundamental form on a minimal surface. Although analogues of the Simons inequality exist for critical points of (\ref{PEF}), the constants degenerate with the ellipticity of $\Phi$, and it is not clear that the same strategy would prove Conjecture \ref{GrowthConj}.



\end{document}